\documentclass[reqno,12pt]{amsart}
\theoremstyle{plain}
\newtheorem{thm}{Theorem}[section]
\newtheorem{lemma}[thm]{Lemma} %%Delete [thm] to re-start numbering
\newtheorem{prop}[thm]{Proposition}

\newtheorem{hyp}[thm]{Hypothesis}

\theoremstyle{remark}

\newtheorem{remark}[thm]{Remark}

\theoremstyle{definition}

\newtheorem{defi}[thm]{Definition}

\usepackage{amscd,amssymb,comment,epic,eepic,euscript,graphics}
\newcount\theTime
\newcount\theHour
\newcount\theMinute
\newcount\theMinuteTens
\newcount\theScratch
\theTime=\number\time
\theHour=\theTime
\divide\theHour by 60
\theScratch=\theHour
\multiply\theScratch by 60
\theMinute=\theTime
\advance\theMinute by -\theScratch
\theMinuteTens=\theMinute
\divide\theMinuteTens by 10
\theScratch=\theMinuteTens
\multiply\theScratch by 10
\advance\theMinute by -\theScratch

\def\today{{\number\day\space
 \ifcase\month\or
  January\or February\or March\or April\or May\or June\or
  July\or August\or September\or October\or November\or December\fi
 \space\number\year}}

\newcommand{\la}{\lambda}
\newcommand{\De}{\Delta}

\newcommand\Rfr{{\mathfrak R}}
\newcommand\Mcal{{\mathcal{M}}} %  needed to be renamed, because was ``\Mc already defined''
\newcommand\Reals{{\mathbf R}}
\newcommand\Complex{{\mathbf C}}
\newcommand\Nats{{\mathbf N}}
\newcommand{\BH}{\mathcal{B}(\mathcal{H})}
\newcommand{\Hcal}{\mathcal{H}}

\newcommand{\st}{\,:\,}
\newcommand{\im}{\text{\rm Im}}
\newcommand{\re}{\text{\rm Re}}
\renewcommand{\i}{\text{\rm i}}

\newcommand{\norm}[1]{\left\Vert#1\right\Vert}

\newcommand\ncs[1]{{\mathcal{L}_#1(\Mcal,\tau)}}
\newcommand{\dd}[3]{\De_{\la_1,\dots,#1}^{(#2)}\left(#3\right)}

\begin{document}

\title[Higher order spectral shift]{Higher order spectral shift, II. Unbounded case.}

\author[Skripka]{Anna Skripka}
\address{Department of Mathematics, Texas A\&M University,
College Station, TX 77843-3368, USA}
\email{askripka@math.tamu.edu}



\date{15 January, 2009}

\begin{abstract}
We construct higher order spectral shift
functions, which represent the remainders of Taylor-type approximations for the
value of a function at a perturbed self-adjoint operator by derivatives of the function at an initial unbounded operator. In the
particular cases of the zero and the first order approximations, the
corresponding spectral shift functions have been constructed by M. G. Krein
\cite{Krein1} and L. S. Koplienko \cite{Kop84}, respectively. The higher order
spectral shift functions obtained in this paper can be expressed recursively
via the lower order ones, in particular, Krein's and Koplienko's spectral shift functions. This extends the recent results of \cite{ds} for bounded operators.
\end{abstract}

\maketitle

\section{Introduction}

Let $\Hcal$ be a separable Hilbert space and $\BH$ the algebra of bounded
linear operators on $\Hcal$. Let $\Mcal$ be a semi-finite von Neumann algebra
acting on $\Hcal$ and $\tau$ a semi-finite normal faithful trace on $\Mcal$. An
initial (possibly) unbounded self-adjoint operator affiliated with $\Mcal$ is
denoted by $H_0$ and its bounded self-adjoint perturbation in $\Mcal$ by $V$.
We continue the study of the recent paper \cite{ds} on the remainder of the
noncommutative Taylor-type approximation
\begin{align}\label{f-la:1}
R_{p,H_0,V}(f)=f(H_0+V)-\sum_{j=0}^{p-1}\frac{1}{j!}\frac{d^j}{dt^j}
\bigg|_{t=0}f(H_0+tV)
\end{align} of the value $f(H_0+V)$, for $f:\Reals\to\Complex$ a bounded function such that the mapping $H\mapsto f(H)$ defined on $H_0+\BH_{sa}$ is $p$ times continuously differentiable in the sense of Fr\'{e}chet (and, hence, in the sense of G\^ateaux). For more details on the mentioned Fr\'{e}chet and G\^ateaux differentiability one can consult \cite{Azamov0,Daletskii,PellerMult} or the introductory section in \cite{ds}.

Let $\mathcal{W}_p$ denote the set of functions $f\in C^p(\Reals)$ such that for
each $j=0,\dots,p$, the derivative $f^{(j)}$ equals the Fourier transform
$\int_\Reals e^{\i t\la}\,d\mu_{f^{(j)}}(\la)$ of a finite Borel measure $\mu_{f^{(j)}}$.
It is known that the functionals $\tau[R_{1,H_0,V}(\cdot)]$ and $\tau[R_{2,H_0,V}(\cdot)]$ are given by absolutely continuous real Borel measures, with the densities $\xi_{H_0+V,H_0}$ and $\eta_{H_0,H_0+V}$ known as Krein's and Koplienko's spectral shift functions, respectively.
More precisely, when $\tau(|V|)<\infty$,
\begin{align}\label{f-la:4}
\tau[R_{1,H_0,V}(f)]=\int_\Reals f'(t)\xi_{H_0+V,H_0}(t)\,dt,
\end{align} for $f\in\mathcal{W}_1$ \cite{Krein1} (see also \cite{Azamov,Carey,Lifshits,PellerKr}), and when $\tau(|V|^2)<\infty$,
\begin{align}\label{f-la:5}
\tau[R_{2,H_0,V}(f)]=\int_\Reals f''(t)\eta_{H_0,H_0+V}(t)\,dt,
\end{align} for $f\in\mathcal{W}_2$ \cite{Kop84} (see also \cite{ds,Neidhardt,PellerKo,convexity}).

It was conjectured in \cite{Kop84} that for $V$ in the Schatten $p$-class, $p\geq 3$, and
$\Mcal=\BH$, there exists a real Borel measure $\nu_p$ with the total variation
bounded by $\frac{\tau(|V|^p)}{p!}$ such that
\begin{align}\label{f-la:7}
\tau[R_{p,H_0,V}(f)]=\int_\Reals f^{(p)}(t)\,d\nu_p(t),\end{align} for bounded
rational functions $f$. Unfortunately, the proof of \eqref{f-la:7} in
\cite{Kop84} has a mistake (see \cite{ds} for details). It was obtained in
\cite[Theorem 5.1]{ds} that \eqref{f-la:7} holds when $V$ is in the
Hilbert-Schmidt class and $\Mcal=\BH$. The total variation of the real Borel
measure $\nu_p$ is bounded by $\frac{\tau(|V|^2)^{p/2}}{p!}$. Analogous results
were also obtained for $\Mcal$ a general semi-finite von Neumann algebra, with
$p=3$, \cite[Theorem 5.2]{ds} and $\Mcal$ a finite von Neumann algebra
and $H_0$ and $V$ free in $\Mcal$, with $p\geq 3$, \cite[Theorem 5.6]{ds}.

Under the assumptions imposed in \cite{ds} in the proof of existence of a finite
measure $\nu_p$ satisfying \eqref{f-la:7}, it was also proved in \cite{ds} that
the measure $\nu_p$ is absolutely continuous when $H_0$ is bounded and its
density $\eta_p$, called the $p$th-order spectral shift function, is
computable recursively by
\begin{align}
\label{f-la:8}
\eta_p(t)=\frac{\tau(V^{p-1})}{(p-1)!}-\nu_{p-1}((-\infty,t))+\int_{\Reals^{p-1}}
spline_{\la_1,\dots,\la_{p-1}}(t)\,dm_{p-1,H_0,V}(\la_1,\dots,\la_{p-1}),
\end{align}
where $spline_{\la_1,\dots,\la_{p-1}}(\cdot)$ is a piecewise polynomial of
degree $p-2$ with breakpoints $\la_1,\dots,\la_{p-1}$ and
$dm_{p-1,H_0,V}(\la_1,\dots,\la_{p-1})$ is a measure on $\Reals^{p-1}$
determined by $p-1$ copies of the spectral measure of $H_0$ intertwined with
$p-1$ copies of the perturbation $V$ (see \eqref{nup} in section \ref{sec:2} for
the precise formula).

Existence of a finite measure $\nu_p$ satisfying \eqref{f-la:7} was proved
first for a bounded operator $H_0$ \cite{ds}. Existence of a finite measure
$\nu_p$ satisfying \eqref{f-la:7} for an unbounded $H_0$ was proved in
\cite{ds} by approximating $H_0$ by bounded operators, applying \eqref{f-la:7}
for bounded operators, and then applying Helly's selection theorem, which gave
no information about the absolutely continuous part of the measure $\nu_p$. In the case
of $p=2$, existence of Koplienko's spectral shift function (that is, the
function $\eta_{H_0,H_0+V}$ satisfying \eqref{f-la:5}) for Hilbert-Schmidt
(non-trace class) perturbations was proved in \cite{Kop84} in the case of
$\Mcal=\BH$ and then extended in \cite{ds} to the case of a general $\Mcal$ by employing the
Birman-Solomyak spectral averaging formula \cite{Birman72}
and the Birman-Solomyak bound
$\big(\tau\big[|f(H_0+V)-f(H_0)|^2\big]\big)^{1/2}\leq
\norm{f'}_\infty\big(\tau\big[|V|^2\big]\big)^{1/2}$,
holding for $f$ a rational function vanishing at infinity (see, e.g.,
\cite{Birman03}). The proof of existence of $\eta_{H_0,H_0+V}$ did not imply
directly that the function $\eta_{H_0,H_0+V}$ is integrable in the case of unbounded $H_0$. Only comparing \eqref{f-la:7}, known to hold for a finite measure $\nu_2$ (whose existence in the case of unbounded $H_0$ is guaranteed by Helly's selection theorem), with \eqref{f-la:5}
implied that $\eta_{H_0,H_0+V}$ is the density of $\nu_2$.

In this paper, we prove (see Theorem \ref{prop:5}) that $d\nu_p(t)=\eta_p(t)dt$
for $H_0$ unbounded, with $\eta_p$ given by \eqref{f-la:8}. The higher
order analogs of the Birman-Solomyak spectral averaging formula exist (see
\cite[Section 10]{ds}), but suitable higher versions of the
Birman-Solomyak bound are not known. Our proof of the
absolute continuity of $\nu_p$ is based on the analysis of the integral in
\eqref{f-la:8} (see section \ref{sec:2}) and the Cauchy (also known as the Borel) transform of the absolutely continuous measure with the density given by \eqref{f-la:8} (see
section \ref{sec:3}), which is not a priori known to be an $L^1$-function. The
obtained equality of the Cauchy transforms of the measures $d\nu_p(t)$ and $\eta_p(t)\,dt$ implies that $d\nu_p(t)=\eta_p(t)dt$. The latter implies (see Remark \ref{rem:to5}) that
$\eta_p$ is integrable and
\begin{align}\label{f-la:9}
\nu_{p-1}(\Reals)=\frac{\tau(V^{p-1})}{(p-1)!}.\end{align} The equality \eqref{f-la:9} was obtained in \cite{Azamov,Krein1} for $p=2$, in \cite{Kop84,PellerKo} for $p=3$, and in \cite{ds} for $p>3$ and $H_0$ bounded.

Let $\ncs{p}$ denote the Schatten
$p$-class of $(\Mcal,\tau)$, that is, $\ncs{p}=\{V\in\Mcal\st\tau(|V|^p)<\infty\}$. The Schatten $p$-class is equipped with the norm
$\norm{\cdot}_{p,\infty}=\norm{\cdot}_p+\norm{\cdot}$, where $\norm{V}_p=\tau(|V|^p)^{1/p}$ and $\norm{\cdot}$ is the operator norm. Throughout the paper, $H_0$ and $V$ denote self-adjoint
operators in $\Mcal$ or affiliated with $\Mcal$; $V$ is mainly taken to be an
element of $\ncs{p}$. Let $\Rfr$ denote the set of rational functions on
$\Reals$ with nonreal poles, $\Rfr_b$ the subset of $\Rfr$ of bounded
functions. The symbol $f_z$ is reserved for the function
$\Reals\ni\la\mapsto\frac{1}{z-\la}$, where $z\in\Complex\setminus\Reals$.

\section{Preliminaries}

\label{sec:1}

In this section we collect facts on divided differences and the Cauchy (Borel) transform of a measure to be used in the proof of the main result.

\subsection{Divided differences}

\begin{defi}\label{prop:dddef}
The divided difference of order $p$ is an operation on functions $f$ of one
(real) variable, which we will usually call $\la$, defined recursively as
follows:
\begin{align*}
&\Delta^{(0)}_{\la_1}(f):=f(\la_1),\\
&\dd{\la_{p+1}}{p}{f}:=
\begin{cases}
\frac{\Delta^{(p-1)}_{\la_1,\dots,\la_{p-1},\la_p}(f)-
\Delta^{(p-1)}_{\la_1,\dots,\la_{p-1},\la_{p+1}}(f)}{\la_p-\la_{p+1}}&
\text{ if } \la_p\neq\la_{p+1}\\[2ex]
\frac{\partial}{\partial t}\big|_{t=\la_p}\Delta^{(p-1)}_{\la_1,\dots,\la_{p-1},t}(f)&
\text{ if } \la_p=\la_{p+1}.
\end{cases}
\end{align*}
\end{defi}

Below we state selected facts on the divided difference.

\begin{prop}
\label{prop:dds}
\begin{enumerate}
\item\label{f-la:sim} $($See \cite[Section 4.7, (a)]{Vore}.$)$\\ $\dd{\la_{p+1}}{p}{f}$
    is symmetric in $\la_1,\la_2,\dots,\la_{p+1}$.

\item\label{f-la:dd2}
$($See \cite[Section 4.7]{Vore}.$)$\\
$\dd{\la_{p+1}}{p}{a_p\la^p+a_{p-1}\la^{p-1}+\dots +a_1\la+a_0}=a_p$, where
$a_0,a_1,\dots,a_p\in\Complex$.

\item\label{f-la:dd3}
$($See \cite[Section 4.7]{Vore}.$)$
For $f$ a sufficiently smooth function,
\[\dd{\la_{p+1}}{p}{f}=\sum_{i\in
\mathcal{I}}\sum_{j=0}^{m(\la_i)-1}c_{ij}(\la_1,\dots,\la_{p+1})
f^{(j)}(\la_i).\] Here $\mathcal{I}$ is
the set of indices $i$ for which $\la_i$ are distinct, $m(\la_i)$
is the multiplicity of $\la_i$, and $c_{ij}(\la_1,\dots,\la_{p+1})\in\Complex$.

\item\label{f-la:dd5} $($See \cite[Section 5.2, (2.3) and (2.6)]{Vore}.$)$

The basic spline with the break points $\la_1,\dots,\la_{p+1}$, where at
least two of the values are distinct, is defined by
\[t\mapsto\begin{cases}
\frac{1}{|\la_2-\la_1|}\chi_{(\min\{\la_1,\la_2\},\max\{\la_1,\la_2\})}(t)&
\text{ if } p=1\\[2ex]
\dd{\la_{p+1}}{p}{(\la-t)^{p-1}_+}& \text{ if } p>1\end{cases}.\] Here the
truncated power is defined by $x_+^k=\begin{cases}x^k &\text{ if }x\geq
0\\0&\text{ if }x<0,\end{cases}$ for $k\in\Nats$.

The basic spline is non-negative, supported in
\[[\min\{\la_1,\dots,\la_{p+1}\},\max\{\la_1,\dots,\la_{p+1}\}]\] and
integrable with the integral equal to $1/p$.

\item\label{f-la:dd4} $($See \cite[Section 5.2, (2.2) and Section 4.7, (c)]{Vore}.$)$

For $f\in
C^p[\min\{\la_1,\dots,\la_{p+1}\},\max\{\la_1,\dots,\la_{p+1}\}]$,
\begin{align}\nonumber&\dd{\la_{p+1}}{p}{f}\\&\quad=\label{f-la:dd*}
\begin{cases}
\frac{1}{(p-1)!}\int_\Reals
f^{(p)}(t)\dd{\la_{p+1}}{p}{(\la-t)^{p-1}_+}\,dt&\;\text{ if
}\;\exists i_1,i_2\text{ such that } \la_{i_1}\neq \la_{i_2}\\[2ex]
\frac{1}{p!}f^{(p)}(\la_1)&\;\text{ if
}\;\la_1=\la_2=\cdots=\la_{p+1}.
\end{cases}\end{align} The first equality in \eqref{f-la:dd*} also holds for $f\in C^{p-1}[a,b]$, with $f^{(p-1)}$ absolutely continuous and $f^{(p)}$ integrable on $[a,b]=[\min\{\la_1,\dots,\la_{p+1}\},\max\{\la_1,\dots,\la_{p+1}\}]$.
\end{enumerate}
\end{prop}

The operation of integration or differentiation can be interchanged with the
operation of taking the divided difference.

\begin{prop}$($See \cite[Lemma 2.4 (iii)]{ds}.$)$\label{prop:dd_der}
Let $D$ be a domain in $\Complex$ and $f$ a function continuously
differentiable sufficiently many times on $D\times\Reals$.
Then for $p\in\Nats$,
\[\frac{\partial}{\partial
z}\left[\dd{\la_{p+1}}{p}{f(z,\la)}\right]=\dd{\la_{p+1}}{p}{\frac{\partial}{\partial
z}f(z,\la)},\] where the divided difference is taken with respect to the
variable $\la$.
\end{prop}

\subsection{The Cauchy transform of a measure}

Let $G_\nu$ denote the Cauchy (also called Borel) transform of a (finite or infinite) measure $\nu$ defined on the Borel subsets of $\Reals$ and satisfying $\int_\Reals\frac{1}{t^2+1}\,d|\nu|(t)<\infty$:
\[G_\nu(z)=\int_\Reals\left(\frac{1}{z-t}+\frac{t}{t^2+1}\right)d\nu(t),\quad \im(z)\neq 0.\]

\begin{prop}$($See, e.g., \cite[Chapter VI, Section 59]{Akhiezer} or \cite[Lemma 7.1]{Kotani}.$)$\label{prop:C1}
If $h$ is a holomorphic function mapping $\Complex_+$ to $\Complex_+$, then it
has the representation \[h(z)=a+bz-G_\nu(z),\] where $a=\re(h(\i))$, $b=\lim_{y\rightarrow\infty}\frac{h(\i y)}{\i y}$,
and $\nu$ is a non-negative measure satisfying\;
$\int_\Reals\frac{1}{t^2+1}\,d\nu(t)<\infty$.
\end{prop}

\begin{remark}\label{rem:C2}
If $\lim_{y\rightarrow\infty}(\i y)h(\i y)\in\Reals$, then $a=\int_\Reals\frac{t}{t^2+1}\,d\nu(t)$ and $b=0$.
\end{remark}

\begin{prop}$($See, e.g., \cite[Chapter VI, Section 59]{Akhiezer} or \cite[Lemma 7.2]{Kotani}.$)$\label{prop:C3}
Given a real measure $\nu$ on $\Reals$ and any $c,d\in\Reals$,
\[\nu((c,d])=-\frac{1}{\pi}\lim_{\delta\rightarrow 0^+}\lim_{\varepsilon\rightarrow 0^+}
\int_{c+\delta}^{d+\delta}\im(G_\nu(t+\i\varepsilon))\,dt.\]
The density of the absolutely continuous part of $\nu$ is given by
\begin{align*}%\label{f-la:Cdens}
\frac{d\nu_{\text\rm ac}(t)}{dt}=-\frac{1}{\pi}\lim_{\varepsilon\rightarrow 0^+}\im(G_\nu(t+\i\varepsilon)),\quad\text{\rm for a.e. }t\in\Reals.
\end{align*}
In particular, if $\im(G_\nu(\cdot))$ is bounded in $\Complex_+$, then the measure $\nu$ is absolutely continuous.
\end{prop}

\begin{lemma}\label{prop:Cparts}
Let $\nu$ be a finite real-valued measure on $\Reals$. Then,
\begin{align*}
G_{\nu}(z)-\int_\Reals\frac{t}{t^2+1}\,d\nu(t)=
\frac{d}{dz}\left(\int_\Reals\left(\frac{1}{z-t}+\frac{t}{t^2+1}\right)
\nu((-\infty,t))\,dt\right),\quad\im(z)\neq 0.\end{align*}
\end{lemma}

\begin{proof}
By performing integration by parts in the representation \[G_\nu(z)-\int_\Reals\frac{t}{t^2+1}\,d\nu(t)=\int_\Reals\frac{1}{z-t}\,d\nu(t),\] we obtain
\begin{align}\nonumber
G_\nu(z)-\int_\Reals\frac{t}{t^2+1}\,d\nu(t)
&=\left(\frac{1}{z-t}\int_{-\infty}^t\,d\nu(\la)\right)\bigg|_{-\infty}^\infty
-\int_\Reals\frac{1}{(z-t)^2}\left(\int_{-\infty}^t\,d\nu(\la)\right)\,dt\\\label{f-la:Cparts1}
&=\frac{d}{dz}\left(\int_\Reals\left(\frac{1}{z-t}+\frac{t}{t^2+1}\right)
\nu((-\infty,t))\,dt\right).
\end{align}
The integral in \eqref{f-la:Cparts1} exists since every bounded function is integrable with the weight $\frac{1}{z-t}+\frac{t}{t^2+1}=\frac{1}{t^2+1}\frac{1+zt}{z-t}$.
\end{proof}

We have the following relation between the remainders of different order.

\begin{prop}$($See \cite[Lemma 3.7]{ds}.$)$\label{prop:R_rec}
Let $H_0=H_0^*$ be an operator in $\Hcal$ and $V=V^*\in\BH$.
%Denote $f_z(t):=\frac{1}{z-t}$, for $z\in\Complex\setminus\Reals$.
Then
\[R_{p+1,H_0,V}\left(f_z\right)=R_{p,H_0,V}\left(f_z\right)
-\big((zI-H_0)^{-1}V\big)^p(zI-H_0)^{-1}.\]
\end{prop}

As a consequence, we obtain a useful relation between the Cauchy transforms
of the measures $\nu_p$ and $\nu_{p+1}$ satisfying \eqref{f-la:7}.

\begin{lemma}\label{prop:Crec}
Let $H_0=H_0^*$ be an operator affiliated with $\Mcal$ and $V=V^*\in\ncs{2}$.
Then for $\nu_p$ satisfying \eqref{f-la:7} for $f=f_z$, with $p\geq 2$,
\begin{align*}
&(-1)^p \tau[R_{p,H_0,V}(f_z)]=G_{\nu_p}^{(p)}(z),\\
&(-1)^{p+1}\tau[R_{p+1,H_0,V}(f_z)]=-G_{\nu_p}^{(p)}(z)+
\frac{(-1)^{p+1}}{p}\frac{d}{dz}\tau\big[\big((zI-H_0)^{-1}V\big)^p\big].\end{align*}
\end{lemma}

\begin{proof} By direct computations, $(-1)^p\tau[R_{p,H_0,V}(f_z)]=G_{\nu_p}^{(p)}(z)$. From Lemma \ref{prop:R_rec} it follows that
\[(-1)^{p+1}\tau[R_{p+1,H_0,V}(f_z)]=-G_{\nu_p}^{(p)}(z)-
(-1)^{p+1}\tau\big[\big((zI-H_0)^{-1}V\big)^p(zI-H_0)^{-1}\big].\] Applying the
representation \[-\tau\big[\big((zI-H_0)^{-1}V\big)^p(zI-H_0)^{-1}\big]=
\frac{d}{dz}\left(\frac{1}{p}\tau\big[\big((zI-H_0)^{-1}V\big)^p\big]\right)\]
completes the proof of the lemma.
\end{proof}

\begin{remark}
The assertion of Lemma \ref{prop:Crec} also holds for $p=1$, provided $V\in\ncs{1}$.
\end{remark}

\section{Properties of splines}

\label{sec:2}

In this section we study properties of piecewise polynomials, which appear in representations for the higher order spectral shift functions.

\begin{lemma}\label{prop:03}
(i) If $\la_1=\dots=\la_p\in\Reals$, with $p\geq 1$, then
\begin{align}\label{f-la:03.-2}
\dd{\la_p}{p-1}{(\la-t)_+^{p-1}}=\chi_{(-\infty,\la_1)}(t).
\end{align}
(ii) If not all $\la_1,\dots,\la_p\in\Reals$ coincide, let $\mathcal{I}$ be
the set of indices $i$ for which $\la_i$ are distinct and let $m(\la_i)$ be
the multiplicity of $\la_i$. Assume that $p\geq 2$ and $M=\max_{i\in\mathcal{I}}m(\la_i)\leq p-1$. Then, $t\mapsto\dd{\la_p}{p-1}{(\la-t)_+^{p-1}}\in C^{p-1-M}(\Reals)$ and
\begin{align}\label{f-la:03.-1}
\dd{\la_p}{p-1}{(\la-t)_+^{p-1}}=(p-1)\int_t^\infty\dd{\la_p}{p-1}{(\la-s)_+^{p-2}}\,ds.
\end{align}
\end{lemma}

\begin{proof}(i) For $t\in\Reals\setminus\{\la_1\}$, the representation \eqref{f-la:03.-2} follows from Proposition \ref{prop:dds} \eqref{f-la:dd4}. For $t=\la_1$ and $p>1$, the claim follows from definition of the divided difference. In the case of $t=\la_1$ and $p=1$, we define $\dd{\la_p}{p-1}{(\la-\la_1)_+^{p-1}}$ to be $0$.

(ii) Now assume that not all $\la_1,\dots,\la_p\in\Reals$ coincide.
The functions $t\mapsto (\la-t)_+^{p-1}$ and $\la\mapsto (\la-t)_+^{p-1}$ are in $C^{p-2}(\Reals)$. By Proposition \ref{prop:dds} \eqref{f-la:dd3},
\begin{align}\label{f-la:03.1}\dd{\la_p}{p-1}{(\la-t)_+^{p-1}}=\sum_{i\in
\mathcal{I}}\sum_{j=0}^{m(\la_i)-1}d_{ij}(\la_1,\dots,\la_p)
(\la_i-t)_+^{p-1-j},\end{align} with $d_{ij}(\la_1,\dots,\la_p)\in\Reals$.
The expression in \eqref{f-la:03.1} is differentiable $p-1-M$ times with respect to $t$ and \[\frac{d^{p-1-M}}{dt^{p-1-M}}\left(\dd{\la_p}{p-1}{(\la-t)_+^{p-1}}\right)=
\dd{\la_p}{p-1}{\frac{\partial^{p-1-M}}{\partial t^{p-1-M}}\left((\la-t)_+^{p-1}\right)}.\]
Finally, Proposition \ref{prop:dds} \eqref{f-la:dd4} implies representation \eqref{f-la:03.-1}.
\end{proof}

\begin{lemma}\label{prop:02}
Let $(\la_1,\dots,\la_p)\in\Reals^p$. Then the function $\dd{\la_p}{p-1}{(\la-t)_+^{p-1}}$ is decreasing; it is equal to $1$ when $t<\min_{1\leq k\leq p}\la_k$ and equal to $0$ when
$t\geq\max_{1\leq k\leq p}\la_k$.
\end{lemma}

\begin{proof}
In the case of $\la_1=\dots=\la_p$, the claim follows from \eqref{f-la:03.-2}. In the case when some of $\la_1,\dots,\la_p$ are distinct, the claim follows from \eqref{f-la:03.-1} and Proposition \ref{prop:dds} \eqref{f-la:dd4}.
\end{proof}

The main assumptions of the following results are collected in the format of a hypothesis.

\begin{hyp}\label{hyp}Let $H_0=H_0^*$ be affiliated with $\Mcal$ and $V=V^*\in\ncs{2}$.
Assume that one of the following three conditions is satisfied:
\begin{enumerate}
\item \label{one} $\Mcal=\BH$, $p\geq 3$,
\item \label{two} $p=3$,
\item \label{three} $\Mcal$ is finite, $p\geq 3$, and $(zI-H_0)^{-1}$ and $V$ are free in $(\Mcal,\tau)$.
\end{enumerate}
\end{hyp}

\begin{prop}$($See \cite[Theorem 5.1 (i), Theorem 5.2 (i), and Theorem 5.6 (i)]{ds}.$)$\label{prop:finite}
Under Hypothesis \ref{hyp} \eqref{one} or \eqref{three}, there exists a unique finite real-valued measure $\nu_p$ such that the trace formula
\begin{equation}\label{f-la:tr}\tau[R_{p,H_0,V}(f)]=\int_\Reals f^{(p)}(t)d\nu_p(t)
\end{equation} holds for $f\in\mathcal{W}_p(\Reals)$. Under Hypothesis \ref{hyp} \eqref{two}, there exists a unique real-valued measure $\nu_p$ such that \eqref{f-la:tr} holds for $f\in C_c^\infty(\Reals)$.
\end{prop}

Under Hypothesis \ref{hyp}, define the function
\begin{align}\label{nup}
\eta_p(t)&=\frac{\tau(V^{p-1})}{(p-1)!}-\nu_{p-1}((-\infty,t))\\\nonumber&\quad-
\frac{1}{(p-1)!}\int_{\Reals^{p-1}}\dd{\la_{p-1}}{p-2}{(\la-t)_+^{p-2}}\,
dm_{p-1,H_0,V}(\la_1,\dots,\la_{p-1}),
\end{align} with $\nu_{p-1}$ a finite measure provided by Proposition \ref{prop:finite} for $p>3$ or by \cite{ds,Kop84} for $p=3$. Here $m_{p-1,H_0,V}$ is the countably additive finite measure on $\Reals^{p-1}$, which uniquely extends the set function
\begin{align*}
A_1\times A_2\times\cdots\times A_p\mapsto
\tau\big[E_{H_0}(A_1)VE_{H_0}(A_2)V\dots VE_{H_0}(A_p)V\big],
\end{align*} where $A_j$, $1\leq j\leq p$, are measurable subsets of $\Reals$ (see \cite[Theorems 4.1 and 4.5]{ds}; also cf. \cite{Birman96,Pavlov}).

\begin{lemma}\label{prop:01}
Let $H_0=H_0^*$ be an operator affiliated with $\Mcal$ and $V=V^*\in\ncs{2}$.
Then $\eta_p$ given by \eqref{nup} is a bounded function satisfying
\[\lim_{t\rightarrow -\infty}\eta_p(t)=0\quad\text{and}\quad
\lim_{t\rightarrow\infty}\eta_p(t)=\frac{\tau(V^{p-1})}{(p-1)!}-\nu_{p-1}(\Reals).\]
\end{lemma}

\begin{proof}
The boundedness of $\eta_p$ follows directly from \eqref{nup}.
The divided difference in \eqref{nup} is symmetric in $\la_1,\dots,\la_{p-1}$ and
\[\overline{m_{p,H_0,V}(d\la_1,d\la_2,\dots,d\la_{p-1})}
=m_{p,H_0,V}(d\la_{p-1},\dots,d\la_2,d\la_1);\] hence, the integral in \eqref{nup} is real-valued (see \cite[proof of Theorem 6.1]{ds} for more details).
We apply Lemma \ref{prop:02} to split the integral in \eqref{nup} into
\begin{align}\nonumber
&\int_{\Reals^{p-1}}\dd{\la_{p-1}}{p-2}{(\la-t)_+^{p-2}}\,
dm_{p-1,H_0,V}(\la_1,\dots,\la_{p-1})\\\label{f-la:02.1}
&\quad=\int_{(t,\infty)^{p-1}}\,dm_{p-1,H_0,V}(\la_1,\dots,\la_{p-1})\\\label{f-la:02.2}
&\quad\quad+\sum\int_{A_{t,k_1}\times\cdots\times A_{t,k_{p-1}}}
\dd{\la_{p-1}}{p-2}{(\la-t)_+^{p-2}}\,dm_{p-1,H_0,V}(\la_1,\dots,\la_{p-1}),
\end{align} where the summation in \eqref{f-la:02.2} goes over all sets $A_{t,k_1}\times\cdots\times A_{t,k_{p-1}}$ such that $A_{t,k_i}\in\{(t,\infty),(-\infty,t]\}$ for every $1\leq i\leq p-1$ and there exist $i_1\neq i_2$ such that $A_{t,k_{i_1}}\neq A_{t,k_{i_2}}$.
The integral in \eqref{f-la:02.1} equals
\[m_{p-1,H_0,V}\left((t,\infty)^{p-1}\right)=
\tau\big[\big(E_{H_0}((t,\infty))V\big)^{p-1}\big],\] which converges to $\tau(V^{p-1})$ as $t\rightarrow-\infty$ and converges to $0$ as $t\rightarrow\infty$ since $E((t,\infty))$ converges to $I$ in the strong operator topology as $t\rightarrow-\infty$ and converges to $\emptyset$ as $t\rightarrow\infty$. We also have
\begin{align*}
&\left|\int_{(t,\infty)\times(-\infty,t]^{p-2}}
\dd{\la_{p-1}}{p-2}{(\la-t)_+^{p-2}}\,dm_{p-1,H_0,V}(\la_1,\dots,\la_{p-1})\right|\\
&\quad\leq|m_{p-1,H_0,V}|\big((t,\infty)\times(-\infty,t]^{p-2}\big),
\end{align*} which is bounded by $|m_{p-1,H_0,V}|\big(\Reals\times(-\infty,t]^{p-2}\big)$ decreasing to $0$ as $t\rightarrow-\infty$ since $\Reals\times(-\infty,t]^{p-2}\searrow\emptyset$. Analogous estimates hold for every set $A_{t,k_1}\times\cdots\times A_{t,k_{p-1}}$ appearing in \eqref{f-la:02.2}. Therefore, the integral in \eqref{f-la:02.2} converges to $0$ as $t\rightarrow-\infty$. Similarly, the integral in \eqref{f-la:02.2} converges to $0$ as $t\rightarrow\infty$. Thus, the integral in \eqref{nup} converges to
$\frac{\tau(V^{p-1})}{(p-1)!}$ as $t\rightarrow-\infty$ and to $0$ as $t\rightarrow\infty$, which along with the asymptotics for \eqref{f-la:02.1} completes the proof of the lemma.
\end{proof}

\section{The main result.}

\label{sec:3}

The theorem below is an extension of the analogous results of \cite{ds} for a bounded operator $H_0$ to the case of an unbounded $H_0$.

\begin{thm}\label{prop:5}
Under Hypothesis \ref{hyp} \eqref{one} or \eqref{three}, the measure $\nu_p$ satisfying \eqref{f-la:tr} for $f\in\mathcal{W}_p$ is absolutely continuous and its density $\eta_p$ is given by \eqref{nup}.
Under Hypothesis \ref{hyp} \eqref{two}, the trace formula
\begin{equation}\label{f-la:treta}
\tau[R_{p,H_0,V}(f)]=\int_\Reals f^{(p)}(t)\eta_p(t)\,dt,
\end{equation} holds for $f\in\Rfr_b$.
\end{thm}

To prove Theorem \ref{prop:5}, we need two auxiliary lemmas.

\begin{lemma}\label{prop:3}
Assume Hypothesis \ref{hyp}. Then,
\begin{align*}
\tau[((zI-H_0)^{-1}V)^p]&=\frac{d^p}{dz^p}\bigg((-1)^p\frac{\tau(V^p)}{(p-1)!}\i\pi
+\frac{(-1)^p}{(p-1)!}\int_\Reals
\bigg(\frac{1}{z-t}+\frac{t}{t^2+1}\bigg)\\
&\quad\quad\quad\quad\int_{\Reals^p}\dd{\la_p}{p-1}{(\la-t)_+^{p-1}}\,
dm_{p,H_0,V}(\la_1,\dots,\la_{p})\,dt\bigg).\end{align*}
\end{lemma}

\begin{proof}As a direct consequence of \cite[Lemma 4.9]{ds}, we have
\begin{align*}
\tau\left[((zI-H_0)^{-1}V)^p\right]=
\int_{\Reals^{p}}
\dd{\la_{p}}{p-1}{\frac{1}{z-\la}}\,dm_{p,H_0,V}(\la_1,\la_2,\dots,\la_p).
\end{align*}
Applying Proposition \ref{prop:dd_der} provides
%(as it was done in the proof of \cite[(40)]{ds}) yields
\begin{align}\label{f-la:3.2}
&\tau\left[((zI-H_0)^{-1}V)^p\right]\\\nonumber
&=\int_{\Reals^{p}}
\dd{\la_{p}}{p-1}{\frac{d^p}{dz^p}\left(\int\dots\int\frac{1}{z-\la}\,dz^p\right)}\,
dm_{p,H_0,V}(\la_1,\la_2,\dots,\la_p)\\
\nonumber
&=\frac{1}{(p-1)!}\int_{\Reals^p}
\frac{d^p}{dz^p}\left(\dd{\la_{p}}{p-1}{(z-\la)^{p-1}\log(z-\la)}+\alpha_{p-1}\right)\,
dm_{p,H_0,V}(\la_1,\la_2,\dots,\la_p),
\end{align} with %a suitable choice of constants of integration on the left hand side and
$\alpha_{p-1}\in\Reals$.

Assume first that not all $\la_1,\dots,\la_p$ coincide.
By Proposition \ref{prop:dds} \eqref{f-la:dd4} and \eqref{f-la:dd5},
\begin{align}\label{f-la:3.3}
&\dd{\la_{p}}{p-1}{(z-\la)^{p-1}\log(z-\la)}+\alpha_{p-1}\\\nonumber
&\quad=(-1)^{p-1}(p-1)\int_\Reals\log(z-t)\dd{\la_{p}}{p-1}{(\la-t)_+^{p-2}}\,dt
+\gamma_{p-1}+\alpha_{p-1},
\end{align} with $\gamma_{p-1}\in\Reals$.
Let
\begin{align}\label{f-la:3.4}
J_{\la_1,\dots,\la_p}(z)=\int_\Reals\log(z-t)\dd{\la_{p}}{p-1}{(\la-t)_+^{p-2}}\,dt.
\end{align}
Since the function $z\mapsto\log(z-t)$ maps $\Complex_+$ to $\Complex_+$ and $\dd{\la_{p}}{p-1}{(\la-t)_+^{p-2}}$ is non-negative, we have that $J_{\la_1,\dots,\la_p}(z)$ maps $\Complex_+$ to $\Complex_+$. By Proposition \ref{prop:dds} \eqref{f-la:dd5}, the function $t\mapsto\dd{\la_{p}}{p-1}{(\la-t)_+^{p-2}}$ is integrable, so
$\lim_{y\rightarrow\infty}\frac{J_{\la_1,\dots,\la_p}(\i y)}{\i y}=0$. By Proposition \ref{prop:C1},
\begin{align}\label{f-la:3.5}
J_{\la_1,\dots,\la_p}(z)=a_{\la_1,\dots,\la_p}-
\int_\Reals\left(\frac{1}{z-t}+\frac{t}{t^2+1}\right)\,d\mu(t),
\end{align}
where $a_{\la_1,\dots,\la_p}=\re\left(J_{\la_1,\dots,\la_p}(\i)\right)$ and $\mu$ is a positive measure satisfying $\int_\Reals\frac{d\mu(t)}{t^2+1}<\infty$.
Due to the integrability of the basic spline (see Proposition \ref{prop:dds} \eqref{f-la:dd5}), the normal boundary value of the imaginary part in \eqref{f-la:3.4} equals
\begin{align*}
&\frac{1}{\pi}\lim_{\varepsilon\rightarrow 0^+}\im\left(J_{\la_1,\dots,\la_p}(t+\i\varepsilon)\right)
=\frac{1}{\pi}\lim_{\varepsilon\rightarrow 0^+}\int_\Reals\im\left(\log(t+\i\varepsilon-s)\right)\dd{\la_{p}}{p-1}{(\la-s)_+^{p-2}}\,ds\\
&\quad=\int_\Reals\chi_{(-\infty,0)}(t-s)\dd{\la_{p}}{p-1}{(\la-s)_+^{p-2}}\,ds
=\int_t^\infty\dd{\la_{p}}{p-1}{(\la-s)_+^{p-2}}\,ds,
\end{align*} which equals $\frac{1}{p-1}\dd{\la_{p}}{p-1}{(\la-t)_+^{p-1}}$ by Lemma \ref{prop:03}. By Proposition \ref{prop:C3}, the measure $\mu$ is absolutely continuous and
\begin{align}\label{f-la:3.6}
\frac{d\mu(t)}{dt}=\frac{1}{p-1}\dd{\la_{p}}{p-1}{(\la-t)_+^{p-1}}.
\end{align} By combining \eqref{f-la:3.3} - \eqref{f-la:3.6}, we obtain
\begin{align*}
&\dd{\la_{p}}{p-1}{(z-\la)^{p-1}\log(z-\la)}+\alpha_{p-1}\\\nonumber
&\quad=\beta_{p-1}-
(-1)^{p-1}\int_\Reals\bigg(\frac{1}{z-t}+\frac{t}{t^2+1}\bigg)
\dd{\la_p}{p-1}{(\la-t)_+^{p-1}}\,dt,
\end{align*} with $\beta_{p-1}\in\Reals$, and hence,
\begin{align}\label{f-la:3.7}
&\frac{d^p}{dz^p}\left(\dd{\la_{p}}{p-1}{(z-\la)^{p-1}\log(z-\la)}+\alpha_{p-1}\right)\\\nonumber
&\quad=\frac{d^p}{dz^p}\left((-1)^p\int_\Reals\bigg(\frac{1}{z-t}+\frac{t}{t^2+1}\bigg)
\dd{\la_p}{p-1}{(\la-t)_+^{p-1}}\,dt\right),
\end{align} where the derivatives are uniformly bounded in $\la_1,\dots,\la_p$ since $\dd{\la_p}{p-1}{(\la-t)_+^{p-1}}$ is bounded (see Lemma \ref{prop:02}). If $\la_1=\dots=\la_p$, then
\[\dd{\la_{p}}{p-1}{(z-\la)^{p-1}\log(z-\la)}+\alpha_{p-1}=
(-1)^{p-1}\log(z-\la_1)+\beta_{p-1},\] with $\beta_{p-1}\in\Reals$, and, by analogous arguments, \eqref{f-la:3.7} follows.
Thus, from \eqref{f-la:3.2} and \eqref{f-la:3.7}, we have
\begin{align*}
&\tau\left[((zI-H_0)^{-1}V)^p\right]\\\nonumber
&=\frac{1}{(p-1)!}\int_{\Reals^p}
\frac{d^p}{dz^p}\left(\dd{\la_{p}}{p-1}{(z-\la)^{p-1}\log(z-\la)}+\alpha_{p-1}\right)\,
dm_{p,H_0,V}(\la_1,\la_2,\dots,\la_p)\\\nonumber
&=\frac{d^p}{dz^p}\bigg(\frac{(-1)^p}{(p-1)!}\int_{\Reals^p}
\int_\Reals\bigg(\frac{1}{z-t}+\frac{t}{t^2+1}\bigg)\\\nonumber
&\quad\quad\quad\quad\quad\quad\quad\quad\quad\dd{\la_p}{p-1}{(\la-t)_+^{p-1}}\,dt\,
dm_{p,H_0,V}(\la_1,\dots,\la_{p})\bigg)\\\nonumber
&=\frac{d^p}{dz^p}\bigg(\frac{(-1)^p}{(p-1)!}\int_\Reals
\bigg(\frac{1}{z-t}+\frac{t}{t^2+1}\bigg)\\\nonumber
&\quad\quad\quad\quad\quad\quad\quad\quad\quad\int_{\Reals^p}\dd{\la_p}{p-1}{(\la-t)_+^{p-1}}\,
dm_{p,H_0,V}(\la_1,\dots,\la_{p})\,dt\bigg),
\end{align*} which completes the proof of Lemma \ref{prop:3} upon adding the $p$th derivative of the constant $(-1)^p\frac{\tau(V^p)}{(p-1)!}\i\pi$.
\end{proof}

\begin{lemma}\label{prop:4}Assume Hypothesis \ref{hyp}. Then,
\[(-1)^p\tau[R_{p,H_0,V}(f_z)]=
\frac{d^p}{dz^p}\left(\int_\Reals\left(\frac{1}{z-t}+
\frac{t}{t^2+1}\right)\eta_p(t)\,dt\right),\]where $\eta_p$ is given by \eqref{nup}.
\end{lemma}

\begin{proof}
Combination of Lemma \ref{prop:Crec}, Lemma \ref{prop:Cparts} applied to the measure $\nu_{p-1}$, which is provided by Proposition \ref{prop:finite}, and Lemma \ref{prop:3} gives
\begin{align}\label{f-la:4.1}\nonumber
&(-1)^p\tau[R_{p,H_0,V}(f_z)]\\\nonumber
&\quad=-G_{\nu_{p-1}}^{(p-1)}(z)+\frac{(-1)^p}{p-1}
\frac{d}{dz}\tau\big[\big((zI-H_0)^{-1}V\big)^{p-1}\big]\\
&\quad=\frac{d^p}{dz^p}\bigg(-\frac{\tau(V^{p-1})}{(p-1)!}\i\pi-
\int_\Reals\left(\frac{1}{z-t}+\frac{t}{t^2+1}\right)
\bigg(\nu_{p-1}((-\infty,t))\\\nonumber
&\quad\quad\quad\quad\quad\quad+
\frac{1}{(p-1)!}\int_{\Reals^{p-1}}\dd{\la_{p-1}}{p-2}{(\la-t)_+^{p-2}}\,
dm_{p-1,H_0,V}(\la_1,\dots,\la_{p-1})\bigg)\,dt\bigg).
\end{align}
Applying the representation $\i\pi=-\int_\Reals\left(\frac{1}{z-t}+\frac{t}{t^2+1}\right)\,dt$ in \eqref{f-la:4.1} completes the proof of the lemma.
\end{proof}

Now we turn to the proof of the main result.

\begin{proof}[Proof of Theorem \ref{prop:5}]
{\it Step 1.} Assume Hypothesis \ref{hyp} \eqref{one} or \eqref{three}.
Applying Lemmas \ref{prop:Crec} and \ref{prop:4} and the equality $\lim_{y\rightarrow\infty}\frac{G_{\nu_p}(\i y)}{\i y}=0$ implies
\begin{align}\label{f-la:5.1}
G_{\nu_p}(z)&=a+\int_\Reals\left(\frac{1}{z-t}+\frac{t}{t^2+1}\right)\eta_p(t)\,dt\\\nonumber
&=\re(a)+\int_\Reals\left(\frac{1}{z-t}+\frac{t}{t^2+1}\right)
\left(\eta_p(t)-\frac{\im(a)}{\pi}\right)\,dt,
\end{align} with $a\in\Complex$. Applying Proposition \ref{prop:C3} in \eqref{f-la:5.1}
implies that
\begin{align}\label{f-la:5.3}
d\nu_p(t)=\left(\eta_p(t)-\frac{\im(a)}{\pi}\right)dt
\end{align} and, in particular, that the measure $\nu_p$ is absolutely continuous. By Lemma \ref{prop:01},
\begin{align*}&\lim_{t\rightarrow-\infty}\left(\eta_p(t)-\frac{\im(a)}{\pi}\right)=
\frac{-\im(a)}{\pi},\\
&\lim_{t\rightarrow\infty}\left(\eta_p(t)-\frac{\im(a)}{\pi}\right)=\frac{-\im(a)}{\pi}+
\frac{\tau(V^{p-1})}{(p-1)!}-\nu_{p-1}(\Reals).\end{align*}
Since $\nu_p$ is finite, $\eta_p-\frac{\im(a)}{\pi}=\frac{d\nu_p}{dt}\in L^1(\Reals)$, and hence, $\im(a)=0$, which along with \eqref{f-la:5.3} completes the proof of the theorem under Hypothesis \ref{hyp} \eqref{one} or \eqref{three}.

{\it Step 2.} Assume Hypothesis \ref{hyp} \eqref{two}. Applying Lemma \ref{prop:4} gives
\begin{align}\label{f-la:th.s2}
\tau[R_{p,H_0,V}(f_z)]=\int_\Reals\frac{p!}{(z-t)^{p+1}}\eta_p(t)\,dt=
\int_\Reals f_z^{(p)}(t)\eta_p(t)\,dt,
\end{align} proving \eqref{f-la:treta} for $f=f_z$. Differentiating \eqref{f-la:th.s2} (similarly to as it was done in \cite{ds}) extends \eqref{f-la:treta} to $f(t)=\frac{1}{(z-t)^k}$, $k\in\Nats$, and, subsequently, to the linear combinations of such functions, which form $\Rfr_b$.
\end{proof}

\begin{remark}\label{rem:to5} Assume Hypothesis \ref{hyp} \eqref{one} or \eqref{three}.
It follows from the proof of Theorem \ref{prop:5} that $\eta_p\in L^1(\Reals)$; it also  follows that $\nu_{p-1}(\Reals)=\frac{\tau(V^{p-1})}{(p-1)!}$, provided $p\geq 3$. When $H_0$ is bounded, $\nu_p(\Reals)=\frac{\tau(V^p)}{p!}$ follows from the trace formula \eqref{f-la:7} applied to $f^{(p)}(t)=t^p$ \cite{ds}.
\end{remark}

\begin{remark}In the case of $p=2$, Theorem \ref{prop:5} holds, provided $V\in\ncs{1}$.
Integrability of $\eta_2$ then follows from the proof of Theorem \ref{prop:5} and finiteness of the measure $\nu_2$.
\end{remark}

\subsection*{Acknowledgment}
The author is thankful to the referee for useful suggestions that led to an improved exposition.

\bibliographystyle{plain}

\end{document}